\documentclass[11pt]{amsart}
\renewcommand{\bar}{\overline}

\newcommand{\pa}{\partial}
\newcommand{\ph}{\varphi}

\usepackage{amsmath,amsthm,amscd,latexsym}

\newcommand{\frk}[1]{{\mathfrak{#1}}}

\title[]{On the Curvature Tensor of the Hodge Metric of 
Moduli Space of Polarized
 Calabi-Yau
Threefolds} 
\author[]{Zhiqin Lu} 
\date{Dec. 14, 1997}
\address[Zhiqin Lu]
{Department of Mathematics\\
Columbia University\\
New York, NY 10027}
\email{lu@math.columbia.edu}

\newtheorem{theorem}{Theorem}[section]
\newtheorem{lemma}{Lemma}[section]
\newtheorem{cor}{Corollary}[section]
\newtheorem{prop}{Proposition}[section]

\newtheorem{definition}{Definition}[section]

\theoremstyle{remark}
\newtheorem{rem}{Remark}[section]
\begin{document}
\maketitle
\tableofcontents

\numberwithin{equation}{section}

\section{Introduction}

This paper is the continuation of the paper~\cite{Lu5} of our study of
the Moduli space of polarized Calabi-Yau threefold.

A polarized Calabi-Yau manifold is a pair $(X,\omega)$ of a compact
algebraic manifold  $X$ with zero first Chern class and a K\"ahler form
$\omega\in H^2(X,Z)$. The form $\omega$ is called a  polarization. Let
$U$ be the universal deformation space of $(X,\omega)$.
$U$ is smooth by a theorem of Tian~\cite{T1}. By~\cite{Y1}, we
may assume that each $X'\in U$ is a K\"ahler-Einstein manifold.
i.e.
the associated K\"ahler metric $(g'_{\alpha\bar\beta})$ is Ricci flat. The
tangent space $T_{X'}{U}$ of ${U}$ at $X'$ can be
identified with $H^1(X',T_{X'})_{\omega}$ where
\[
H^1(X',T_{X'})_{\omega}=\{\phi\in H^1(X',T_{X'})| \phi\lrcorner\omega=0\}
\]
The Weil-Petersson metric
$G_{PW}$ on ${U}$ is defined by
\[  
G_{WP}(\phi,\psi)=\int_{X'}
{g'}^{\alpha\bar\beta}g'_{\gamma\bar\delta}\phi^\gamma_{\bar\beta}
\bar{\psi^\delta_{\bar\alpha}}
dV_{g'}
\]
where $\phi=\phi^\gamma_{\bar\beta}\frac{\partial}{\partial z^\gamma}
d\bar{z}^\beta$,
$\psi=\psi^\delta_{\bar\alpha}\frac{\partial}{\partial z^\delta}
d\bar{z}^\alpha\in H^1(X',T_{X'})_\omega$,
and $g'=g'_{\alpha\bar\beta}
dz^\alpha d\bar{z}^\beta$ is the K\"ahler-Einstein metric on $X$
associated to  the polarization $\omega$.

Let $n=\dim U$. As showed in ~\cite{Lu5}, we defined the
Hodge metric $\omega_H$ by
\[
\omega_H=(n+3)\omega_{WP}+Ric (\omega_{WP})
\]
where $\omega_{WP}$ is the K\"ahler form of the Weil-Petersson metric. 

The main result of ~\cite{Lu5} is the following

\begin{theorem}\label{01}
Let $\omega_H=(n+3)\omega_{WP}+Ric(\omega_{WP})$. Then
\begin{enumerate}
\item $\omega_H$ is a K\"ahler metric on $U$;
\item The holomorphic bisectional curvature of $\omega_H$ is nonpositive.
Furthermore,
Let $\alpha=((\sqrt{n}+1)^2+1)^{-1}>0$. Then the  Ricci curvature $Ric
(\omega_H)\leq
-\alpha \omega_H$ and the holomorphic sectional curvature is also less
than or
equal to $-\alpha$.
\item If $Ric(\omega_H)$ is bounded, then the Riemannian sectional
curvature of $\omega_H$ is also bounded.
\end{enumerate}
\end{theorem}

The main result of this paper builds on the above theorem:

\begin{theorem}
\label{thm11}
Let $X$ be a Calabi-Yau threefold. Let $\ph_1,\cdots,\ph_n$
be the orthonormal harmonic basis of $H^1(X,T_{X})_{\omega}$.
Then there is a  constant $C$,
depending only on $n$, such that 
the $L^\infty$ norm
 of the sectional curvature $|R|$ satisfies
\[
|R|\leq C\sum_{i=1}^{n}||\ph_i||^8_{L^4}
\]
\end{theorem}

\begin{rem}
The crucial part of this theorem is that the curvature has an upper bound
which only depends on the $L^4$ norm of the harmonic basis, rather
than depends on the derivative of the harmonic basis. 
Upper bound of the sectional curvature of the Hodge metric is very
important in the compactification of the moduli space
(cf. ~\cite{Lu7}).
\end{rem}

In order to prove the theorem, we  need to estimate the covariant
derivative of the Yukawa coupling with respect to the Weil-Petersson
metric. As a by-product, we proved the following theorem (for
definitions, see \S 2): 

\begin{theorem}\label{thm12}
Let $F=(F_{ijk})$ be the Yukawa coupling. Let
\[
F_{ijk,l}=\pa_l F_{ijk} -\Gamma^m_{il} F_{mjk}
-\Gamma^m_{jl} F_{imk}
-\Gamma^m_{kl} F_{ijm}
+2 K_l F_{ijk}
\]
Then $F_{ijk,l}=F_{ijl,k}$.
\end{theorem}

\begin{rem}
The moduli space of a Calabi-Yau threefold is a {\it projective} special 
K\"ahler manifold in the sense of D. Freed~\cite{DSF}. In~\cite{Lu8}, the
Yakuwa coupling of 
special K\"ahler manifolds is discussed.
\end{rem}

The motivation behind this paper and the paper~\cite{Lu5} is
that we want to  give a differential geometric proof of the theorem of
Viehweg~\cite{V2}
in the case of
the moduli space of Calabi-Yau threefolds. 
Viehweg's theorem states that moduli spaces of polarized algebraic
varieties are quasi-projective.
The boundedness of the
curvature of the Hodge metric is very important because of the work of
Mok~\cite{Mok1}, Mok-Zhong~\cite{MZ} and Yeung~\cite{Yeung1}.
By their theorems, if  a complete K\"ahler manifold of finite volume has
negative Ricci curvature and bounded sectional curvature, then it must be
quasi-projective. 
In the case of the moduli space of Calabi-Yau threefolds, the Ricci
curvature is negative (Theorem~\ref{01}), and the condition on the
boundedness of the sectional curvatures can be weakened, thus it is very 
important to get various upper bound estimates of the sectional 
curvatures.

In the last section, we give an extra restriction
on the limit of Hodge structures for a one dimensional
degeneration of a family of Calabi-Yau threefolds.

\smallskip

\smallskip

{\bf Acknowledgment.}
The author thanks Professor Tian for his constant encouragement 
and discussion during the preparation of this paper. He also thanks C. L.
Wong for a lot of useful conversations.

\section{The Covariant Derivatives of the Yukawa Coupling}
Suppose $\pi: {\frk  X}\rightarrow U$ is the total space over
the (local) universal deformation space $U$ of a Calabi-Yau
threefold $X$. Thus
for any point $X'\in U$, $\pi^{-1}(X')$ is a Calabi-Yau threefold.
The Hodge bundle $\underline{F}^3=\pi_*\omega_{{\frk X}/U}$ is the
holomorphic line bundle over $U$ where $\omega_{{\frk
X}/U}$ is the relative canonical bundle of ${\frk  X}$.

There is a natural Hermitian metric on $\underline{F}^3$ defined by the
Ricci flat metric on each fiber of $\pi$.
Such a metric can be written out explicitly as follows: since for any
$X'\in U$, $\pi^{-1}(X')$ is differmorphic to $\pi^{-1}(0)=X$, there is
a natural identification of $H^3(X',C)\rightarrow H^3(X,C)$. Suppose
$\ph,\psi\in
H^3(X,C)$. Define the cup product
\[
Q(\ph,\psi)=-\int_X \ph\wedge\psi
\]

Let $\Omega$ be a local holomorphic section of $\underline{F}^3$.
Thus for each $X'\in U$, $\Omega$ at $X'$ is a holomorphic $(3,0)$ form
on $H^{3,0}(X')\subset H^3(X',C)$ and under the identification
$H^3(X',C)\rightarrow H^3(X,C)$, $\Omega(X')\in H^3(X,C)$.

The Hermitian metric on $\underline{F}^3$ is defined by setting 
$||\Omega||^2=\sqrt{-1} Q(\Omega,\bar\Omega)$.

The technical heart of this paper is to compute the covariant derivative
of the Yukawa coupling with respect to the Weil-Petersson metric
and the Hermitian metric on the bundle $\underline{F}^3$. Recall
that by definition(see\cite{BG}, for example), the Yukawa coupling is the
(local) section $F$ 
of the bundle
$Sym^3((R_*^1(T_{{\frk X}/U}))^*)\otimes (\underline{F}^3)^{\otimes 2}$ 
over $U$ such that for any $\ph_1,\ph_2,\ph_3\in H^1(X',T_{X'})$ and
$\Omega\in H^{3,0}(X')$, 
\[
F(\ph_1,\ph_2,\ph_3)=
\int_{X'} (\ph_1\wedge\ph_2\wedge\ph_3\lrcorner\Omega)\wedge\Omega
\]
Here $T_{{\frk X}/U}$ is the relative tangent sheaf of
${\frk X}\rightarrow U$.

The basic property of the Yukawa coupling is that it is a holomorphic
section. In fact, Let $t^1,\cdots, t^n$ be the local holomorphic
coordinate system of $U$. Let $\Omega$ be a local nonzero section of the
holomorphic
bundle $\underline{F}^3$, 
We have
\[
F_{ijk}=F(\rho(\frac{\pa}{\pa t^i}),\rho(\frac{\pa}{\pa
t^j}),\rho(\frac{\pa}{\pa t^k}))
=Q(\Omega,\pa_i\pa_j\pa_k\Omega),\qquad 1\leq i,j,k\leq n
\]
where $\rho : T_XU\rightarrow H^1(X,T_X)$ is the Kodaira-Spencer
map.

Let $\Gamma_{ij}^k$ be the Christoffel symbols of the K\"ahler metric
$g_{i\bar j}$ and let $K_l=-\pa_l\log ||\Omega||^2$ be the connection of
the Hermitian bundle
$\underline{F}^3$
with respect to the local section $\Omega$. 
We make the following definition:

\begin{definition}
For $1\leq i,j,k\leq n$, the covariant derivative  
of $F_{ijk}$ is defined as
\[
F_{ijk,l}=\pa_l F_{ijk}-\Gamma^s_{li} F_{sjk}
-\Gamma^s_{lj} F_{isk}
-\Gamma^s_{lk} F_{ijs}
+2 F_{ijk} K_l
\]
\end{definition}

\vspace{0.2in}

In this section, we are going  to compute $F_{ijk,l}$ at a point
$X\in U$
in terms of the information of the fixed Calabi-Yau threefold $X$.

We use the method developed by Siu\cite{Siu2}, Nannicini\cite{Nnn} and
Schumacher\cite{Sc1}. 
By K\"ahler geometry, there is a holomorphic
coordinate $(t^1,\cdots, t^n)$ of $U$ such that at $X$, $\Gamma^k_{ij}=0$,
$1\leq i,j,k\leq n$. Furthermore, if the local section $\Omega$ of
$\underline{F}^3$ is carefully chosen, then $K_l=0$, $1\leq l\leq n$ at
$X$.

Consider the Kodaira-Spencer map $\rho : T_{X'}U\rightarrow
H^1(X',T_{X'})$. Let $\ph_j=\rho(\frac{\pa}{\pa t^j}), 1\leq j\leq n$.
Suppose  $\ph_j$'s are harmonic $T_{X'}$-valued $(0,1)$ forms. These
$\ph_j$'s can
be realized
by the canonical lift in the sense of Siu~\cite{Siu2} (See also
Nannicini~\cite{Nnn} and Schumacher~\cite{Sc1}): 
suppose $(z^1,z^2,z^3)$
is the holomorphic coordinate on $X$. Then for each $\frac{\pa}{\pa t^j}$,
there is a vector $v_j$, called the canonical lift of $\frac{\pa}{\pa
t^j}$, on ${\frk X}$ which locally can be represented as
$v_j=\frac{\pa}{\pa t^j}+v^\alpha_j\frac{\pa}{\pa z^\alpha}$ such that
$\ph_j=\bar\pa v_j^\alpha \frac{\pa}{\pa z^\alpha}$
is a harmonic $T_{X'}$-valued (0,1) form.

It should be noted that $v_j$ is a vector field on ${\frk X}$ but neither
is $\frac{\pa}{\pa t^j}$ nor $v_j^\alpha \frac{\pa}{\pa z^\alpha}$ alone.
In fact, the component $\frac{\pa}{\pa t^j}$ in the expression 
$v_j=\frac{\pa}{\pa t^j}+v^\alpha_j\frac{\pa}{\pa z^\alpha}$
is  different from   it is as the vector field on $U$. 
It is also easy to check that the real part of the vector field $v_j$
defines differmorphisms between the fibers. Using these differmorphisms, 
tensor fields of the nearby fibers can be identified as tensor fields
on $X$. 
By Nannicini~\cite{Nnn} or Siu~\cite{Siu2}, the Lie derivative
$L_{v_l}\cdot$ is defined as the usual $\frac{\pa}{\pa t^l}$
after pulling back via the differmorphisms.

Now let's analyze the conditions $K_l=0$ and $\Gamma^i_{jk}=0$ for $1\leq
i,j,k,l\leq n$ at $X$. We have

\begin{prop}\label{prop21}
We use the notations as above. In particular, suppose
$(\frac{\pa}{\pa t^1},\cdots, \frac{\pa}{\pa t^n})$ and $\Omega$ are 
chosen such that
$\Gamma_{jk}^i=0$ and $K_l=0$ at $X\in U$. Then
we have
\begin{enumerate}
\item $(L_{v_l}\Omega)^{3,0}=0$;
\item $(L_{v_k}\ph_i)$ is a $\bar\pa^*$-boundary.
\end{enumerate}
\end{prop}

{\bf Proof:}
The key point is to identify the derivatives with respect to 
cohomological classes and the derivatives with respect to  forms. 
Suppose
for fixed $l$, $v_l=\tau_1+\sqrt{-1}\tau_2$
where $\tau_1$ and $\tau_2$ are real vector fields.
 Let $\sigma_1(s)$ and
$\sigma_2(s)$ be the flows defined by $\tau_1$ and $\tau_2$,
respectively.
Consider the 3-forms $\sigma_i(s)^*\bar\Omega$, $i=1,2$. Suppose
\[
\sigma_i(s)^*\bar\Omega=p_i(s)+dq_i(s)
\]
be the Hodge decomposition of $\sigma_i(s)^*\bar\Omega$ in $H^3(X,C)$.
Then we have
\[
0=\pa_l\bar\Omega=\frac{d}{ds}|_{s=0}(p_1(s)+\sqrt{-1}p_2(s))
\]
by the definition of $\pa_l\bar\Omega$. This is equivalent to 
\[
\frac{d}{ds}|_{s=0} (\sigma_1(s)^*\bar\Omega
+\sqrt{-1}\sigma_2(s)^*\bar\Omega)
-d\frac{d}{ds}|_{s=0} (q_1(s)+\sqrt{-1}q_2(s))=0
\]
Or in other word
\[
L_{v_l}\bar\Omega-d\sigma=0
\]
for $\sigma=\frac{d}{ds}|_{s=0}(q_1+\sqrt{-1}q_2)$.

Using this, we have
\[
\int_X\Omega\wedge L_{v_l}\bar\Omega=0
\]

On the other hand, $K_l=0$ implies
\[
0=\pa_l Q(\Omega,\bar\Omega)
=\int_X L_{v_l}\Omega\wedge\bar\Omega
+\int_X\Omega\wedge L_{v_l}\bar\Omega
\]

So the first part of the proposition follows from the following:

{\bf Claim.\,} $\bar\pa L_{v_l}\Omega=0$.

{\bf Proof of the Claim:} This follows from a straightforward 
computation. Let  $\Omega$ be represented as
\[
\Omega=a dz^1\wedge dz^2\wedge d z^3
\]
where the functions $a, z^1,z^2,z^3$ are holomorphic on each fiber
and have parameter $t$. Suppose $\rho(\frac{\pa}{\pa t})=\ph
=\ph^\alpha_{\bar\beta}\pa_\alpha d\bar z^{\beta}$ is a
harmonic $T_X$-valued $(0,1)$ form.
We have
\begin{align}\label{a}
\begin{split}
& \bar\pa (\frac{\pa}{\pa t} (dz^1\wedge dz^2\wedge dz^3))\\
&=\bar\pa d\frac{\pa z^1}{\pa t}\wedge dz^2\wedge dz^3
+\bar\pa (dz^1\wedge d\frac{\pa z^2}{\pa t}\wedge dz^3)
+\bar\pa (dz^1\wedge dz^2\wedge d\frac{\pa z^3}{\pa t})\\
&=-\pa\bar\pa \frac{\pa z^1}{\pa t}\wedge dz^2\wedge dz^3
+ dz^1\wedge \pa\bar\pa\frac{\pa z^2}{\pa t}\wedge dz^3
-dz^1\wedge dz^2\wedge \pa\bar\pa\frac{\pa z^3}{\pa t}\\
&=(-\pa_1\ph^1_{\bar\beta}-\pa_2\ph^2_{\bar\beta}
-\pa_3\ph^3_{\bar\beta})
dz^1\wedge dz^2\wedge dz^3\wedge d\bar{z^\beta}
\end{split}
\end{align}

By the harmonicity of $\ph$,
we have
\begin{equation}\label{b}
\pa_\alpha\ph^\alpha_{\bar\beta}
+\Gamma^\alpha_{\alpha\gamma}\ph^\gamma_{\bar\beta}
=0
\end{equation}
where the notation $\Gamma^\alpha_{\beta\gamma}$ is the connection of $X$
which is
different from the connection $\Gamma^i_{jk}$ on the universal deformation
space $U$.
 
From the theory of deformation of complex structures, we know
that $\bar\pa-t\ph$ defines the $\bar\pa$-operator on the nearby fibers. 
Thus we have
\begin{equation}\label{c}
\bar\pa\frac{\pa a}{\pa t} 
=\ph^\alpha_{\bar\beta}\pa_\alpha a d\bar{z^\beta}
\end{equation}

Using Equation~(\ref{a}), (\ref{b}) and (\ref{c}), we have

\begin{align*}
& \bar\pa\frac{\pa}{\pa t}\Omega
=\bar\pa \frac{\pa a}{\pa t} dz^1\wedge dz^2\wedge dz^3
+a\bar\pa\frac{\pa}{\pa t}(dz^1\wedge dz^2\wedge dz^3)\\
&=(\ph^\alpha_{\bar\beta}\pa_\alpha a-a\Gamma^\alpha_{\alpha\gamma}
\ph^\gamma_{\bar\beta}) d\bar z^\beta\wedge dz^1\wedge dz^2\wedge dz^3
\end{align*}

So $\bar\pa\frac{\pa}{\pa t}\Omega=0$ follows from the fact
\[
a\Gamma^\alpha_{\alpha\gamma}=a\pa_\gamma\log |a|^2
=a\pa_\gamma\log a=\pa_\gamma a
\]
and the claim is proved.

The second 
part of the proposition  is implied in ~\cite{Nnn}. We prove it for the
sake
of completeness. 
We assume at $X$, $(z^1,z^2,z^3)$ are normal coordinates.
By definition,
\[
L_{v_k}\ph_i=
(\pa_k(\ph_i)^\alpha_{\bar\beta}
-(\ph_i)^\gamma_{\bar\beta}\pa_\gamma v_k^\alpha)\pa_\alpha d\bar z^\beta
\]

Thus
\begin{equation}\label{x}
\bar\pa^* L_{v_k}\ph_i=(\pa_\beta\pa_k (\ph_i)^\alpha_{\bar\beta}
-(\ph_i)^\gamma_{\bar\beta}\pa_\beta\pa_\gamma v_k^\alpha)\pa_\alpha
\end{equation}

We are going to prove $\bar\pa^* L_{v_k}\ph_i=0$. By the harmonicity of
$\ph_i$, we have
\[
g^{\beta_1\bar\beta}(\ph_i)^\alpha_{\bar\beta,\beta_1}=0
\]
or
\[
g^{\beta_1\bar\beta}(\pa_{\beta_1}(\ph_i)^\alpha_{\bar\beta}
+\Gamma^\alpha_{\beta_1\beta_2}(\ph_i)^{\beta_2}_{\bar\beta})=0
\]

Taking derivative with respect to $\pa_k$ gives
\begin{equation}\label{y}
g^{\beta_1\bar\beta}\pa_{\pa_k\beta_1}(\ph_i)^\alpha_{\bar\beta}+
\pa_k g^{\beta_1\bar\beta}\pa_{\beta_1}(\ph_i)^\alpha_{\bar\beta}
+\pa_k (g^{\beta_1\bar\beta}\Gamma^\alpha_{\beta_1\beta_2})
(\ph_i)^{\beta_2}_{\bar\beta}=0
\end{equation}

Since $L_{v_k}\omega=0$ (See Nannicini~\cite{Nnn}, for example), we have
\[
0=L_{v_k} g_{\beta\bar\beta_1} dz^\beta\wedge d\bar z^{\beta_1}
=(\pa_k g_{\beta\bar\beta_1}+\pa_\beta v_k^{\beta_1}) dz^\beta\wedge d
\bar z^{\beta_1}
\]

So we have
\begin{equation}\label{d}
\pa_k g_{\beta\bar\beta_1}=-\pa_\beta v_k^{\beta_1}
\end{equation}

We also have
\begin{equation}\label{e}
\pa_k(g^{\beta_1\bar\beta}\Gamma^\alpha_{\beta_1\beta_2})
=\pa_k\Gamma^\alpha_{\beta\beta_2}
=\pa_k\Gamma_{\beta\bar\alpha\beta_2}
=\pa_k\pa_{\beta_2} g_{\beta\bar\alpha}
=-\pa_{\beta_2}\pa_\beta v_k^\alpha
\end{equation}

In addition, we have
\begin{equation}\label{f}
\pa_k\pa_{\beta_1}(\ph_i)^\alpha_{\bar\beta}
=\pa_{\beta_1}\pa_k (\ph_i)^\alpha_{\bar\beta}
-\pa_{\beta_1} v_k^\gamma\pa_\gamma (\ph_i)^\alpha_{\bar\beta}
\end{equation}

Using Equation~(\ref{d}), (\ref{e}) and (\ref{f}), from 
Equation~(\ref{x}) and (\ref{y}), we get $\bar\pa^* L_{v_k}\ph_i=0$.

By~\cite{Nnn}, we see that
\[
0=\Gamma_{i\bar jk}=\int_X <L_{v_k}\ph_i, \ph_j>
\]

So the harmonic part of $L_{v_k}\ph_i$ is zero. Thus $L_{v_k}\ph_i$
is a $\bar\pa^*$-boundary by the Hodge decomposition.

\qed

The condition that $(\frac{\pa}{\pa t^1},\cdots, \frac{\pa}{\pa t^n})$
is an orthonormal basis implies that
 $\ph_1,\cdots,\ph_n\in
H^1(X,T_X)$ is a set of orthnormal basis of harmonic $T_X$-valued
forms.
Let $\Omega$ be the local  nonzero section of $\underline{F}^3$. We
make the following definition:

\begin{definition}\label{d21}
\[
a_{jk}=(\ph_j\wedge \ph_k)^\#=\ph_j\wedge\ph_k\lrcorner\Omega
\]
is an (1,2) form for $1\leq j,k\leq n$.
Here 
\[
\ph_j\wedge\ph_k
=(\ph_j)^\alpha_{\bar\beta}(\ph_k)^\gamma_{\bar\delta}
\frac{\pa}{\pa z^\alpha}\wedge\frac{\pa}{\pa z^\gamma}
\otimes d\bar z^\beta\wedge d\bar z^\delta
\]
If $\Omega=a dz^1\wedge dz^2\wedge dz^3$. Then $a_{jk}$ can be represented
as
\begin{equation}\label{m}
a_{jk}=a(\ph_j)^\alpha_{\bar\beta}(\ph_k)^\gamma_{\bar\delta}
 sgn(\xi,\alpha,\gamma) dz^\xi\wedge d\bar z^\beta\wedge d\bar z^\delta
\end{equation}
\end{definition}

\begin{lemma}\label{lem21} For $1\leq j,k\leq n$,
\[
\pa^*a_{jk}=0
\]
\end{lemma}

{\bf Proof:} Since $X$ is a K\"ahler manifold, we have
$\pa^*=\sqrt{-1}[\Lambda,\bar\pa]$. First we have
$\bar\pa a_{jk}=0$. Next, suppose $(z^1,z^2,z^3)$ is a normal coordinate
system.
Then
\begin{align*}
& \Lambda a_{jk}=
\Lambda (a (\ph_j)^\alpha_{\bar\beta}(\ph_k)^\gamma_{\bar\delta}
sgn(\alpha,\gamma,\xi) dz^\xi\wedge d\bar{z}^{\beta}\wedge
d\bar{z}^\delta)\\
&=a(\ph_j)^\alpha_{\bar\beta}(\ph_k)^\gamma_{\bar\delta}
sgn (\alpha,\gamma,\beta) d\bar{z}^\delta
-a(\ph_j)^\alpha_{\bar\beta}(\ph_k)^\gamma_{\bar\delta}
sgn (\alpha,\gamma,\delta) d\bar{z}^\beta
\end{align*}

However, the fact that $\ph_j\lrcorner\omega=0$ and
$\ph_k\lrcorner\omega=0$ implies
$(\ph_j)^\alpha_{\bar\beta}=(\ph_j)^\beta_{\bar\alpha}$, and
$(\ph_k)^\gamma_{\bar\delta}=(\ph_k)^\delta_{\bar\gamma}$.
So $\Lambda a_{jk}=0$ and the lemma is proved.

\qed

\begin{definition}
The  Hodge *-operator on $X$ is defined as
\[
*: A^{p,q}\rightarrow A^{n-p,n-q},\qquad
(\ph, \psi) dV=\ph\wedge *\psi
\]
for $\ph,\psi\in A^{p,q}(X)$.
\end{definition}

\begin{lemma}\label{lem22}
For $1\leq j,k\leq n$
\[
*a_{jk}=\bar{a_{jk}}
\]
\end{lemma}

{\bf Proof:}
By Equation~(\ref{m}),
\[
a_{jk}=a(\ph_j)^\alpha_{\bar\beta}(\ph_k)^\gamma_{\bar\delta}
sgn (\alpha, \gamma, \xi) dz^\xi\wedge d\bar{z}^\beta\wedge
d\bar{z}^\delta
\]
We have
\begin{align*}
& *a_{jk}=
\sum_{m<n}\bar{a(\ph_j)^\alpha_{\bar\beta}(\ph_k)^\gamma_{\bar\delta}}
sgn (\alpha,\gamma,\xi) sgn(\xi, m, n) sgn(\beta,\delta,\eta)
d\bar z^m\wedge d\bar z^n\wedge dz^\eta\\
& =\bar{a(\ph_j)^\alpha_{\bar\beta}(\ph_k)^\gamma_{\bar\delta}}
d\bar z^\alpha\wedge d\bar z^\gamma sgn (\beta, \delta, \eta)
d\bar{z}^\eta
=\bar{a_{jk}}
\end{align*}
Here we again use the fact $\ph_j\lrcorner\omega=0$ and
$\ph_k\lrcorner\omega=0$.

\qed

\begin{lemma}\label{lem23} 
For $1\leq j,k\leq n$,
\[
\bar\pa (L_{v_j} \ph_k\lrcorner\Omega)=\pa a_{jk}
\]
\end{lemma}

{\bf Proof:} In Nannicini~\cite{Nnn}, it is proved that
\[
\bar\pa L_{v_j}\ph_k= D_1^* (\ph_j\wedge\ph_k)
\]
where
\[
D_1^*(\ph_j\wedge\ph_k)
=\pa_\alpha ((\ph_j)^\alpha_{\bar\beta} (\ph_k)^\gamma_{\bar\delta})
\pa_\gamma d\bar{z}^\beta\wedge d \bar z^\delta
\]
The lemma follows from Equation~(\ref{m}).

\vspace{0.1in}

Now we are going to prove the main theorem of this section.

\begin{theorem}\label{thm21}
Let $G$ be the Green's operator on differential forms of $X$. 
Suppose $\ph_1,\cdots,\ph_n$ are the orthnormal basis of
$H^1(X,T_X)$. Then
\[
F_{ijk,l}=\int_X (G\pa a_{li},\bar{\pa a_{jk}}) dV
+\int_X (G\pa a_{lj},\bar{\pa a_{ik}}) dV
+\int_X (G\pa a_{lk},\bar{\pa a_{ij}}) dV
\]
where $a_{ij}$'s are defined in Definition~\ref{d21}.
\end{theorem}

It should be noted that the notation of the inner product is defined as
$(a,a)=||a||^2$. So $(a,\bar a)$ is {\it not} the norm of $a$.

\vspace{0.2in}

{\bf Proof:} 
Since we choose the local coordinate $(t_1,\cdots,t_n)$ and the local
section $\Omega$ such that
$\Gamma^i_{jk}=K_l=0$, we have
$F_{ijk,l}=\pa_l F_{ijk}$. 
By Proposition~\ref{prop21}, 
\begin{align*}
&\pa_l F_{ijk}
=\pa_l\int_X(\ph_i\wedge\ph_j\wedge\ph_k)\lrcorner\Omega)\wedge\Omega
=\int_X (L_{v_l}\ph_i\wedge\ph_j\wedge\ph_k\lrcorner\Omega)\wedge\Omega\\
&+\int_X (\ph_i\wedge L_{v_l}\ph_j\wedge\ph_k\lrcorner\Omega)\wedge\Omega
+\int_X (\ph_i\wedge\ph_j\wedge L_{v_l}\ph_k\lrcorner\Omega)\wedge\Omega
\end{align*}

Thus we need only  to prove that
\[
\int_X ((L_{v_l}\ph_i)\wedge\ph_j\wedge\ph_k\lrcorner\Omega)\wedge\Omega
=\int_X (G\pa a_{li},\bar{\pa a_{jk}}) dV
\]

Let $L_{v_l}\ph_i=b_{li}$. Recall that $\Omega=adz^1\wedge dz^2\wedge
dz^3$. Then
\begin{equation}\label{aa}
(b_{li}\lrcorner\Omega)
=(-1)^\alpha a(b_{li})^\alpha_{\bar\alpha_1}
dz^1\wedge\cdots \hat{dz^\alpha}\cdots \wedge dz^3\wedge d\bar
z^{\alpha_1}
\end{equation}

By Equation~(\ref{m}) we see that
\begin{align*}
&(b_{li}\lrcorner\Omega\wedge a_{jk})
=(-1)^\alpha a^2(b_{li})^\alpha_{\bar\alpha_1}
(\ph_j)^\beta_{\bar\beta_1} (\ph_k)^\gamma_{\bar\gamma_1}
dz^1\wedge\cdots \hat{dz^\alpha}\cdots \wedge dz^3
\wedge d\bar z^{\alpha_1}\\
&\qquad\qquad sgn (\beta,\gamma,\xi) dz^\xi\wedge d\bar
z^{\beta_1}\wedge d\bar
z^{\gamma_1}\\
&=a^2(b_{li})^\alpha_{\bar\alpha_1}
(\ph_j)^\beta_{\bar\beta_1}(\ph_k)^\gamma_{\bar\gamma_1}
sgn(\alpha,\beta,\gamma)
sgn(\alpha_1,\beta_1,\gamma_1)\\
&\qquad\qquad dz^1\wedge dz^2\wedge dz^3\wedge d\bar z^1\wedge d\bar z^2
\wedge d\bar z^3\\
&=-(b_{li}\wedge\ph_j\wedge\ph_k\lrcorner\Omega)\wedge\Omega
\end{align*}

By Lemma~\ref{lem22} and the above equation, we have
\begin{align}\label{aaa}
\begin{split}
& (b_{li}\wedge\ph_j\wedge\ph_k\lrcorner\Omega)\wedge\Omega
=-(b_{li}\lrcorner\Omega\wedge a_{jk})\\
&=(b_{li}\lrcorner\Omega\wedge *\bar{a_{jk}})
=(b_{li}\lrcorner\Omega,\bar{a_{jk}}) dV
\end{split}
\end{align}

By Lemma~\ref{lem23} and Proposition~\ref{prop21}, we have
\begin{align*}
&\Box (L_{v_l}\ph_i\lrcorner\Omega)
=\bar\pa^* \bar\pa (L_{v_l}\ph_i\lrcorner\Omega)
+\bar\pa\bar\pa^* (L_{v_l}\ph_i\lrcorner\Omega)\\
&=\bar\pa^*\pa a_{li}+\bar\pa ((\bar\pa^* L_{v_l}\ph_i)\lrcorner\Omega)
=\bar\pa^*\pa a_{li}
\end{align*}

Since $L_{v_l}\ph_i$ is a $\bar\pa^*$-boundary (Proposition~\ref{prop21}), we
know
\begin{equation}\label{bbb}
L_{v_l}\ph_i\lrcorner\Omega=G\bar\pa^*\pa a_{li}
\end{equation}
where $G$ is the Green operator of the Laplacian $\Box$.
Thus by Equation~(\ref{aaa}) and (\ref{bbb}),
\begin{align*}
&\int_X ((L_{v_l}\ph_i)\wedge\ph_j\wedge\ph_k\lrcorner\Omega)
\wedge\Omega
=\int_X(b_{li}\lrcorner\Omega,\bar{a_{jk}}) dV\\
&=\int_X (G\bar\pa^*\pa a_{li},\bar{a_{jk}})
=\int_X (G\pa a_{li},\bar{\pa a_{jk}}) dV
\end{align*}

\qed

\begin{cor}[Theorem~\ref{thm12}] For $1\leq i,j,k,l\leq n$,
\[
F_{ijk,l}=F_{ijl,k}
\]
\end{cor}
\qed

\section{The Estimates}
In this section, we give an  upper bound of the curvature tensor of the
Hodge metric. We use the same notations as in the previous section.

Suppose $(z^1,\cdots,z^n)$ is the normal holomorphic coordinate 
system at $p\in U$ with
respect to the Weil-Petersson metric
$\omega_{WP}=\sqrt{-1} g_{i\bar j} dz^i\wedge d\bar z^j$.

 We further assume that
$Q(\Omega,\bar\Omega)=1$. In ~\cite{Lu5} we
have proved
that

\begin{theorem}
If $(z^1,\cdots,z^n)$ is the normal coordinate system of $\omega_{WP}$,
then the curvature tensor $R_{i\bar jk\bar l}$ of $\omega_H
=\sqrt{-1}h_{i\bar j} dz^i\wedge d\bar z^j$ at $p$ is 
\[
{R}_{i\bar{j}k\bar{l}}=A_{i\bar{j}k\bar{l}}
+B_{i\bar{j}k\bar{l}}
\]
where
\begin{align*}
&
A_{i\bar{j}k\bar{l}}=2\delta_{ij}\delta_{kl}
+2\delta_{il}\delta_{kj}-4\sum_sF_{iks}\bar{F_{jls}}
+2\sum_{mnpq}F_{qkm}\bar{F_{plm}}F_{inp}\bar{F_{jnq}}\\
&
B_{i\bar{j}k\bar{l}}=
\sum_{rs}F_{irs,k}\bar{F_{jrs,l}}
-\sum_{mn}(\sum_{rs}F_{irs,k}\bar{F_{mrs}})
\bar{(\sum_{rs}F_{jrs,l}\bar{F_{nrs}})}h^{n\bar m}
\end{align*}
Here 
$F_{ijk}$ is the Yukawa coupling and $F_{ijk,l}$ is its covariant
derivative with respect to the Weil-Petersson metric and the connection
on $\underline{F^3}$.
\end{theorem}

It is proved in~\cite{Lu5}  that in order to bound the curvature
tensor, we need
only  to bound the scalar curvature. 
By definition, the scalar curvature $\rho$
is
\[
\rho=
-h^{i\bar j}h^{k\bar l} R_{i\bar jk\bar l}=
-h^{i\bar j}h^{k\bar l} A_{i\bar jk\bar l}
-h^{i\bar j}h^{k\bar l} B_{i\bar jk\bar l}
\]

\begin{lemma}\label{lem31}
Suppose that the dimension of the universal deformation space $U$ is $n$. Then
\[
|h^{i\bar j}h^{k\bar l} A_{i\bar jk\bar l}|
\leq 3n^6
\]
\end{lemma}

{\bf Proof:}
Under the local coordinate $(z^1,\cdots,z^n)$, 
we proved in~\cite{Lu5} that
\[
h_{i\bar j}=2\delta_{i\bar j}+\sum_{mn}F_{imn}\bar{F_{jmn}}
\]
Suppose further that 
\[
\sum_{mn}F_{imn}\bar{F_{jmn}}=\lambda_i\delta_{ij}
\]
Then
\[
h_{i\bar j}=(2+\lambda_i)\delta_{ij}
\]

In particular, for fixed $i,m,n$
\[
|F_{imn}|\leq\sqrt{\lambda_i}
\]

So 
\begin{align*}
&\sum_{ijkl} h^{i\bar j}h^{k\bar l} \sum_{mnpq}
F_{qkm}\bar{F_{plm}}F_{inp}\bar{F_{jnq}}\\
&=\sum_{ikmnpq}\frac{1}{2+\lambda_i}\cdot\frac{1}{2+\lambda_k}
F_{qkm}\bar{F_{pkm}}F_{inp}\bar{F_{inq}}\\
&\leq\sum_{ik} 
\frac{1}{2+\lambda_i}\cdot\frac{1}{2+\lambda_k}\lambda_i\lambda_k n^4
\end{align*}

In~\cite{Lu5}, we have proved that $h^{i\bar j}h^{k\bar l} A_{i\bar jk\bar l}\geq
0$.
Thus
\[
|h^{i\bar j}h^{k\bar l} A_{i\bar jk\bar l}|
\leq
\sum_{ik}\frac{4}{(2+\lambda_i)(2+\lambda_k)}
+2\sum_{ik}
\frac{1}{2+\lambda_i}\cdot\frac{1}{2+\lambda_k}\lambda_i\lambda_k n^4
\leq 3n^6
\]

\qed

Now we consider $h^{i\bar j}h^{k\bar l} B_{i\bar jk\bar l}$. It is easy to
see that
\begin{equation}\label{mn}
0\leq h^{i\bar j}h^{k\bar l} B_{i\bar jk\bar l}\leq
\sum_{ijklrs} h^{i\bar j}h^{k\bar l}
F_{irs,k}\bar{F_{jrs,l}}
\leq\sum_{ijkl}|F_{ijk,l}|^2
\end{equation}

\begin{lemma}\label{lem32} 
Using the notations in Theorem ~\ref{thm21} ,
we have
\[
|F_{ijk,l}|
\leq
3(||\ph_i||^4_{L^4}+||\ph_j||^4_{L^4}+||\ph_k||^4_{L^4}+||\ph_l||^4_{L^4})
\]
\end{lemma}

{\bf Proof:}
Since the Green operator is a positive operator, we have
\[
|\int_X (G\pa a_{li},\bar{\pa a_{jk}})|
\leq
\sqrt{\int_X (G\pa a_{li},\pa a_{li})}
\sqrt{\int_X (G\pa a_{jk},\pa a_{jk})}
\]

However, for fixed $l,i$, by Lemma ~\ref{lem21}, we have

\[
G\pa^*\pa a_{li}=G\Box a_{li}
=a_{li}-H(a_{li})
\]
where $H(a_{li})$ is the harmonic part in the Hodge decomposition 
of $a_{li}$. Thus

\[
|\int_X (G\pa a_{li},\pa a_{li})|
\leq
|\int_X (G\pa^*\pa a_{li},a_{li})|
\leq
\int_X |a_{li}|^2
\leq
\int_X |\ph_l|^4+|\ph_i|^4
\]

So by Theorem~\ref{thm21},
\[
|F_{ijk,l}|
\leq
3(||\ph_i||^4_{L^4}+||\ph_j||^4_{L^4}+||\ph_k||^4_{L^4}+||\ph_l||^4_{L^4})
\]

\qed

\begin{theorem}[Theorem~\ref{thm11}]
The scalar curvature $\rho$ of the Hodge metric  satisfies
\[
0<-\rho\leq 3n^6+144n^3\sum_i||\ph_i||^8_{L^4}
\]
\end{theorem}

{\bf Proof:}
$\rho<0$ follows from ~\cite{Lu5}.
The upper bound is from Lemma~\ref{lem31} and Lemma~\ref{lem32}.
\qed

\section{A Remark on the Theorem of C-L Wang}
In his paper~\cite{Wong}, C-L Wong gave a necessary and sufficient
condition for the
Weil-Petersson metric to be incomplete 
for a family of Calabi-Yau manifolds
over a punctured disk. The main
theorem of him is (for the precise definitions and notations, see
Wong~\cite{Wong}):

\begin{theorem}[C-L Wong~\cite{Wong}]
Let $\Delta^*$ be the parameter space of a family of Calabi-Yau
manifolds. Let $F_\infty^n$ be the limit of $F^n$
in the sense of Hodge theory and $N$
is the associated 
nilpotent operator. Then the necessary and sufficient condition for the
Weil-Petersson metric to be incomplete is that $NF_\infty^n=0$.
\end{theorem}

In this section, we are going to prove, even if the Weil-Petersson metric 
is
complete, we still have some restrictions on $F^n_\infty$ for $n=3$.

The classical Weil-Petersson metric is defined by giving a natural
Hermitian metric
on $H^1(X,T_X)$ induced by the Ricci flat K\"ahler metric
for each Calabi-Yau manifold. However, by the theorem of Tian~\cite{T1},
We can look at the Weil-Petersson metric in a different way.

Recall that the Hodge bundle $\underline{F^n}$ over the classifying space
$D$ is the tautological bundle of the filtration
\[
0\subset F^n\subset F^{n-1}\subset\cdots\subset F^1\subset H
\]
The natural Hermitian metric on $\underline{F^n}$ is the polarization
$Q$. Suppose $\omega$ is the curvature form of the Hermitian metric $Q$,
then $\omega$ is an closed (1,1) form of $D$. Suppose $M$ is a horizontal
slice of $D$ (see Griffiths~\cite{Gr}, for example), then $\omega$
restricts to a semi-positive form on $M$. However, if $M$ is  the
universal deformation space of a Calabi-Yau manifold, then by
Tian's theorem~\cite{T1}, $\omega|_M$ must be positive definite and is the
Weil-Petersson metric.

Thus there are some restrictions for a horizontal slice on which the
$\omega$ is positive definite. 
 The following theorem gives one of the
restrictions on  the limiting Hodge structure.

\begin{theorem} We use the notations in the above theorem and
in~\cite{Wong}.
If $n=3$, then
\[
Q(F_\infty^3, N^3F_\infty^3-3N^2F_\infty^3-2NF_\infty^3)=0
\]
\end{theorem}

{\bf Proof:} Let
\[
\Omega=e^{\frac{1}{2\pi\sqrt{-1}}\log zN} A(z)
\]
where $A(z)$ is a vector valued holomorphic function of $z\in\Delta^*$,
the punctured unit disk. Let
\[
F_{zzz}= (\Omega,\pa_z\pa_z\pa_z\Omega)
\]

It is easy to check  that
\[
\lim_{z\rightarrow 0} z^3 F_{zzz}=
Q(F_\infty^3, N^3F_\infty^3-3N^2F_\infty^3-2NF_\infty^3)
\]

So we need only to prove that
\[
\lim_{z\rightarrow 0} z^3 F_{zzz}=0
\]

Let $p\in\Delta^*$. Then since $p$ represents a Calabi-Yau threefold, we
have a map $f$ from a neighborhood of $p$ in $\Delta^*$ to the universal
deformation space $U$. Suppose in local coordinates, the map $f$ is
$z\mapsto (z^1,\cdots,z^n)$. Let $Z^i=\frac{\pa z^i}{\pa z}$. Then from
~\cite{Lu5}, we see that the Hodge metric on $\Delta^*$ can be written as
\[
h=h_{i\bar j} Z^i\bar{Z}^j
=(2 g_{i\bar j}+g^{m\bar n} g^{p\bar q} F_{imp}\bar{F_{jnq}})
Z^i\bar{Z}^j
\]
where $g_{i\bar j}$ and $h_{i\bar j}$ are the Weil-Petersson
metric and the Hodge metric, respectively. Since
$g_{i\bar j}\leq h_{i\bar j}$, we have
\[
h\geq h^{m\bar n} h^{p\bar q} F_{imp}\bar{F_{jnq}}
Z^i\bar{Z}^j
\]

By the Cauchy inequality, we see that
\[
(h^{m\bar n} h^{p\bar q} F_{imp}\bar{F_{jnq}}
Z^i\bar{Z}^j) h^2\geq
|F_{ijk}Z^iZ^jZ^k|^2
\]
So we have
\[
h^3\geq |F_{ijk}Z^iZ^jZ^k|^2=|F_{zzz}|^2
\]

In ~\cite{Lu5}, it is proved that the curvature of $h$ is 
negative away from zero. So the Schwartz lemma gives,

\[
h\leq\frac{C}{r^2(\log \frac 1r)^2}
\]

Then
\[
r^6|F_{zzz}|^2\leq C\frac{r^6}{r^6 (\log \frac 1r)^6}\rightarrow 0
\]
The theorem is proved.

\qed

\bibliographystyle{acm}
\bibliography{bib}

\end{document}